\documentclass[12pt]{article}
\usepackage[final]{epsfig}
\usepackage{graphics}
\usepackage{amsmath}
\usepackage{amsfonts}
\usepackage{latexsym}
\usepackage{amssymb}
\usepackage{graphicx}
\usepackage{url}
\usepackage{epstopdf}

\newtheorem{lemma}{Lemma}[section]
\newtheorem{proposition}[lemma]{Proposition}

\newtheorem{example}[lemma]{Example}
\newtheorem{theorem}{Theorem}

\newtheorem{conjecture}[lemma]{Conjecture}
\newtheorem{problem}[lemma]{Problem}

\newcommand{\proofend}{$\Box$\bigskip}

\newcommand{\R}{{\mathbb R}}

\newcommand{\HH}{{\mathbb H}}
\newcommand{\Sph}{{\mathbb S}}

\def\proof{\paragraph{Proof.}}
\newcommand{\Vol}{\mathrm{Vol}}

\begin{document}

\title{Remarks on the the circumcenter of mass}

\author{Serge Tabachnikov\footnote{
Department of Mathematics,
Pennsylvania State University,
University Park, PA 16802,
and ICERM, Brown University, Box 1995, Providence, RI 02912;
tabachni@math.psu.edu} \and Emmanuel Tsukerman\footnote{Department of Mathematics, University of California, Berkeley, CA 94720-3840; e.tsukerman@berkeley.edu}
}

\date{}
\maketitle

%\begin{abstract}
%\end{abstract}

\section{Introduction} \label{intro}

Given a homogeneous polygonal lamina $P$, one way to find its center of mass is as follows: triangulate $P$, assign to each triangle its centroid, taken with the weight equal to the area of the triangle, and find the center of mass of the resulting system of point masses. That the resulting point, $CM(P)$, does not depend on the triangulation, is a consequence of the Archimedes Lemma: {\it if an object is divided into smaller objects, then the center of mass of the compound object is the weighted average of the centers of mass of the parts, with the weights equal to the respective areas}. 

\begin{figure}[hbtp]
\centering
\includegraphics[height=2in]{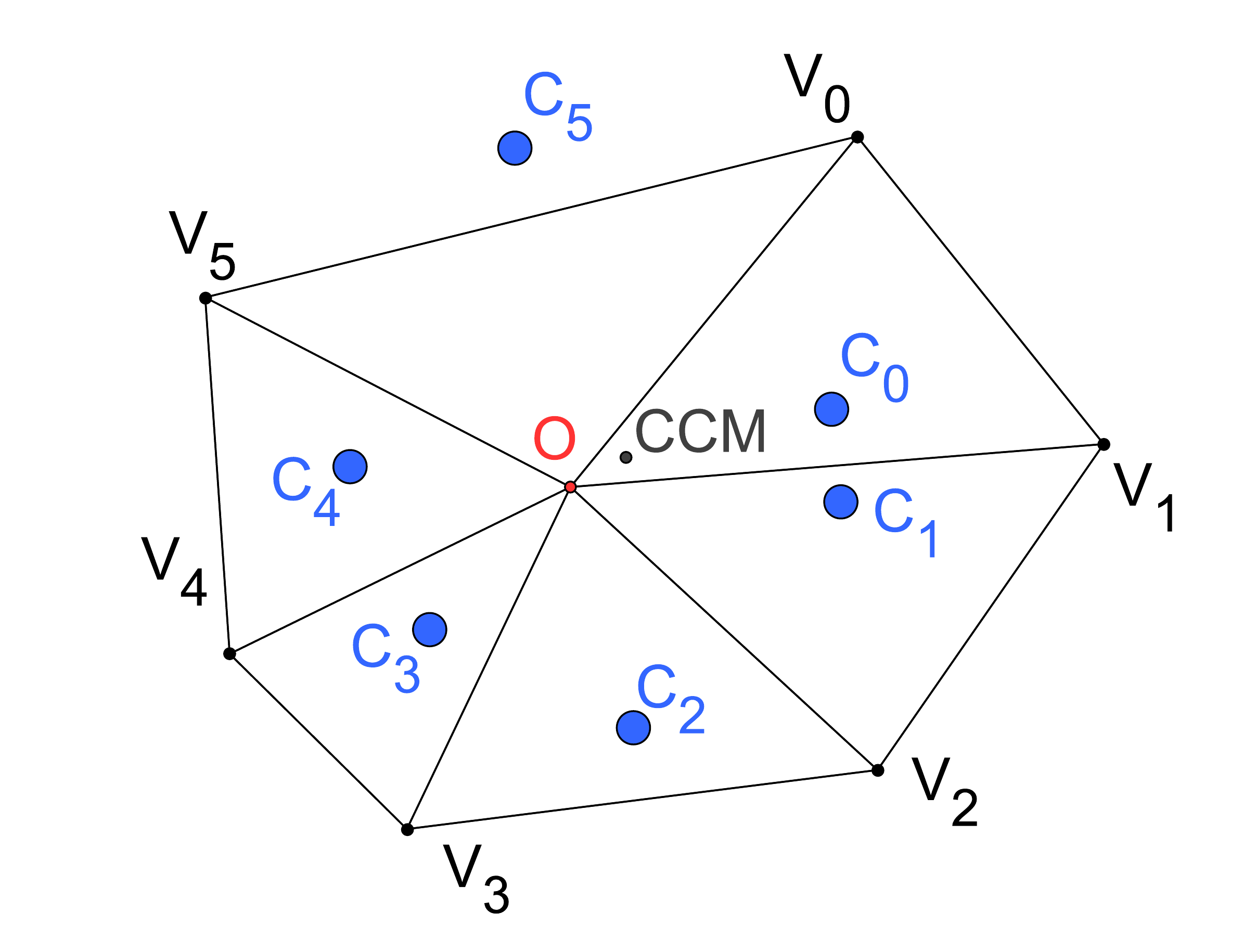}
\caption{Circumcenter of mass}
\label{Defn}
\end{figure}

Replace, in the above construction, the centroids of the triangles by their circumcenters. The resulting weighted average is called the {\it circumcenter of mass} of the polygon $P$, denoted by $CCM(P)$. This point is well defined, that is, does not depend on the triangulation (assuming that degenerate triangles are avoided), see Figure \ref{Defn}.

This construction is mentioned in the 19th century book \cite{La}, where it is attributed to the Italian algebraic geometer G. Bellavitis. We learned about this reference from B. Gr\"unbaum who, together with G. C. Shephard, studied this construction in the early 1990s \cite{Gr}. Independently, and at about the same time, the circumcenter of mass was rediscovered by V. Adler \cite{Ad1,Ad2} as an integral of a discrete dynamical system called  recutting of polygons.

The explicit formulas are as follows. Let the coordinates of the vertices of the polygon $P$, taken in the cyclic order, be $(x_i,y_i),\ i=1,\ldots,n$. Then $CCM(P)=$
$$
\frac{1}{4 A(P)}\left(\sum_{i=0}^{n-1}y_{i}(x_{i-1}^{2}+y_{i-1}^{2}-x_{i+1}^{2}-y_{i+1}^{2}),\sum_{i=0}^{n-1}-x_{i}(x_{i-1}^{2}+y_{i-1}^{2}-x_{i+1}^{2}-y_{i+1}^{2})\right),
$$
where $A(P)$ is the signed area of $P$. For comparison,  $CM(P)=$
$$
\frac{1}{6 A(P)} \left(\sum_{i=0}^{n-1}(x_{i}+x_{i+1})(x_{i}y_{i+1}-x_{i+1}y_{i}),\sum_{i=0}^{n-1}(y_{i}+y_{i+1})(x_{i}y_{i+1}-x_{i+1}y_{i})\right).
$$

The construction of the  circumcenter of mass extends to higher dimensions, and to the elliptic and hyperbolic geometries. We studied it in \cite{TT2} in relation with the so-called discrete bicycle transformation \cite{TT1}. See also the  paper by A. Akopyan \cite{Ak}.

The construction in $\R^n$ is similar. Given  a simplicial polytope $P$, consider its non-degenerate triangulation. Assign the circumcenter $CC(\Delta_i)$ to each simplex  $\Delta_i$ of the triangulation, and take the center of mass of these points with weights equal to the oriented volumes of the respective simplices: 
\begin{equation} \label{totsum}
CCM(P)=\frac{1}{\Vol(P)}\sum_i \Vol(\Delta_i)\ CC(\Delta_i).
\end{equation}
The result does not depend on the triangulation. 

The explicit formula is as follows. Let $F=(V_1,\ldots,V_n)$ be a face of $P$, where $V_i$ are vectors in $\R^n$. Let $A(F)$ be the $n\times n$ matrix made of vectors $V_i$, and let $A_i(F)$ be obtained form $A(F)$ by replacing $i$th row with $(|V_i|^2,\ldots,|V_n|^2)$. Then the $i$th component of the circumcenter of mass is given by
$$
CCM(P)_i=\frac{1}{2 (n!) {\rm Vol}(P)} \sum_{F\subset \partial P} \det A_i(F).
$$

One can take affine combinations $t CM + (1-t) CCM,\ t\in\R$, resulting in a line,  called the {\it generalized Euler line} of the polytope $P$ (for a triangle, the Euler line is the line through the centroid and the circumcenter; it passes through the orthocenter as well).

In this note we are interested in the uniqueness of this construction. 

Suppose that to every non-degenerate simplex $\Delta \subset \R^n$ a `center' $C(\Delta) \in \R^n$ is assigned so that the following assumptions hold:
\begin{enumerate}
\item The map $\Delta \mapsto C(\Delta)$ commutes with similarities (both orientation-preserving and orientation-reversing);
\item The map $\Delta \mapsto C(\Delta)$ is invariant under the permutations of the vertices of the simplex $\Delta$;
\item The map $\varphi: \Delta \mapsto \Vol(\Delta) C(\Delta)$ is polynomial in the coordinates of the vertices of the simplex $\Delta$;
%\item Archimedes Lemma: if a simplex $\Delta$ is triangulated by choosing a point $O$, as described above, then an analog of the equation (\ref{Arch}) holds. 
\end{enumerate}

\begin{theorem} \label{main}
Under these assumptions, $C(\Delta)$ is an affine combination of the center of mass and the circumcenter: 
$$C(\Delta) = t CM(\Delta) + (1-t) CC(\Delta),$$
where the constant $t\in\R$ depends on the map $\Delta \mapsto C(\Delta)$ (and does not depend on the simplex $\Delta$).
\end{theorem}

\section{Basic determinants} \label{basic}

Let $x_1,\ldots,x_n$ be Cartesian coordinates in $\R^n$. Let $\Delta =(V_0,\ldots,V_n)$ be a simplex, and let  $V_j=(x_1^j,\ldots,x_n^j),\ j=0,\ldots,n$, be the coordinates of its vertices (where $j$ is a superscript, not an exponent).
Let 
\begin{align*}
V=\left|\begin{array}{ccccccccccccc}
x_1^0 & x_1^1 & \cdots & x_1^n \\
x_2^0 & x_2^1 & \cdots & x_2^n \\
\vdots & \vdots & \ddots & \vdots \\ 
x_n^0 & x_n^1 & \cdots & x_n^n \\
1 & 1 & \cdots & 1
\end{array}\right|, 
\end{align*}
a multiple of the oriented volume of $\Delta$, and 
\begin{align*}
X_{i,jk}=\left|\begin{array}{ccccccccccccc}
x_1^0 & x_1^1 & \cdots & x_1^n \\
x_2^0 & x_2^1 & \cdots & x_2^n \\
\vdots & \vdots & \ddots & \vdots \\ 
x_{i-1}^0 & x_{i-1}^1 & \cdots & x_{i-1}^n \\
x_j^0 x_k^0 & x_j^1 x_k^1 & \cdots & x_j^n x_k^n \\
x_{i+1}^0 & x_{i+1}^1 & \cdots & x_{i+1}^n \\
\vdots & \vdots & \ddots & \vdots \\
x_n^0 & x_n^1 & \cdots & x_n^n \\
1 & 1 & \cdots & 1
\end{array}\right|=\text{Skew}(x_1^0 x_2^1 \cdots \widehat{x_i^{i-1}} \cdots x_n^{n-1} x_j^{i-1} x_k^{i-1}), \quad 1 \leq i,j,k \leq n.
\end{align*}
where Skew is skew-symmetrization over superscripts.
Evidently, $X_{i,jk}=X_{i,kj}$, and the number of such polynomials equals $n^2(n+1)/2$. Both determinants, $V$ and  $X_{i,jk}$, are skew-symmetric under permutations of the vertices of the simplex.

\begin{lemma} \label{basis}
The polynomials $X_{i,jk}$ constitute a linear basis of the space ${\cal S}$ of homogeneous polynomials of degree $n+1$ in the variables $x_1^0,x_2^0,\ldots,x_n^n$, skew-symmetric under permutations of the superscripts.
\end{lemma}

\proof
Since 
\[
X_{i,jk}=\text{Skew}(x_1^0 x_2^1 \cdots \widehat{x_i^{i-1}} \cdots x_n^{n-1} x_j^{i-1} x_k^{i-1}),
\]
 a monomial $x_1^0 x_2^1 \cdots \widehat{x_i^{i-1}} \cdots x_n^{n-1} x_j^{i-1} x_k^{i-1}$ determines $X_{i,jk}$. We will show that: 
\begin{enumerate} 
 \item There exists no nonidentity permutation acting on superscripts that maps this monomial to itself. 
 \item These monomials lie in different orbits under this action. 
\end{enumerate} 
 The former shows that this monomial does not cancel in the expression of $X_{i,jk}$. The latter will then imply that different monomials give rise to different determinants. 
 
 Suppose that there exists a permutation $\sigma$ such that
\[
x_1^0 x_2^1 \cdots \widehat{x_i^{i-1}} \cdots x_n^{n-1} x_j^{i-1} x_k^{i-1}=x_1^{\sigma(0)} x_2^{\sigma(1)} \cdots \widehat{x_{i'}^{\sigma(i'-1)}} \cdots x_n^{\sigma(n-1)} x_{j'}^{\sigma(i'-1)} x_{k'}^{\sigma(i'-1)}.
\]
The superscripts $i-1$ and $\sigma(i'-1)$ are the unique ones which occur twice. Therefore $\sigma(i'-1)=i-1$. The corresponding subscripts are $j,k$ and $j',k'$, so that $\{j,k\}=\{j',k'\}$. Dividing both sides by $x_j^{i-1} x_k^{i-1}$, we get
  \[
x_1^0 x_2^1 \cdots \widehat{x_i^{i-1}} \cdots x_n^{n-1} =x_1^{\sigma(0)} x_2^{\sigma(1)} \cdots \widehat{x_{i'}^{\sigma(i'-1)}} \cdots x_n^{\sigma(n-1)},
\]
so that $i=i'$ and $\sigma$ is the identity.

Let $f \in {\cal S}$. Then $f$ is equal to its skew-symmetrization. Write $f$ in its monomial basis:
$$
f(x_1^0,\ldots,x_n^n)=\sum c_{\beta}^{\alpha} x_{\beta}^{\alpha}, \quad |\alpha|=|\beta|=n+1.
$$
Consider the skew-symmetrization of a monomial $x_{\beta}^{\alpha}$. If some number appears in $\alpha$ with multiplicity $3$ or greater, then $\alpha$ must be missing some two distinct numbers $i,j \in \{0,1,\ldots,n\}$. Each permutation $\sigma$ has a counterpart $\sigma (i \, j)$ of opposite sign which maps $x_\beta^{\alpha}$ to the same monomial. Therefore the skew-symmetrization of $x_\beta^{\alpha}$ in this case is zero. 

Now suppose that $\alpha$ contains $n+1$ different elements of $\{0,1,\ldots,n\}$. Since the entries of $\beta$ are elements of $\{1,2,\ldots,n\}$, there exist some $\beta_i$ and $\beta_j$ with $i \neq j$ such that $\beta_i=\beta_j$. Each permutation $\sigma$ has a counterpart $\sigma (\alpha_i \, \alpha_j)$ of opposite sign which maps $x_\beta^{\alpha}$ to the same element.    

It follows that the only monomials appearing in $f$ are those for which $\alpha$ is a permutation of $(0,1,\ldots,\widehat{i-1},\ldots,n-1,i-1,i-1)$. Assume without loss of generality that $\alpha$ is of this form. Let $\alpha=(\alpha_1,\alpha_2,\ldots,\alpha_n,\alpha_{n+1})$ and $\beta=(\beta_1,\beta_2,\ldots,\beta_n,\beta_{n+1})$. Suppose that $\beta_{i_1}=\beta_{i_2}=\ldots=\beta_{i_k}$. For the monomial to not vanish under skew-symmetrization, the corresponding multiset $\alpha_{i_1},\alpha_{i_2},\ldots,\alpha_{i_k}$ cannot be invariant under any transpositions. Knowing the structure of $\alpha$, we see that this implies that $k=1,2$ or $3$. If $k=2$, then $(\alpha_{i_1},\alpha_{i_2})=(i-1,i-1)$. If $k=3$ then  $(\alpha_{i_1},\alpha_{i_2},\alpha_{i_3})=(i-1,i-1,q)$, with $q \neq i-1$. This proves the claim.
\proofend

Consider the map $\varphi: \Delta \mapsto \Vol(\Delta) C(\Delta)$, and let $(y_1,\ldots,y_n)$ be its components. Our assumption 3 implies that each $y_\ell,\ \ell=1,\ldots,n$, is a polynomial in the variables $x_1^0,x_2^0,\ldots,x_n^n$. The assumption 1, applied to scaling, implies that these polynomials are homogeneous of degree $n+1$, and the assumption 2 that they are skew-symmetric under permutations of the superscripts. Lemma \ref{basis} implies that
\begin{equation*}
y_\ell = \sum A_{i,jk}^\ell X_{i,jk},\ \ell=1,\ldots,n,
\end{equation*}
where the coefficients $A_{i,jk}^\ell$ satisfy $A_{i,jk}^\ell= A_{i,kj}^\ell$. We always assume that summation is over repeated indices.

\begin{example} \label{twoknown}
{\rm The center of mass and the circumcenter of mass correspond to the functions
$$
y_\ell = \frac{1}{n+1}\sum X_{i,i\ell}\ \ {\rm and}\ \ y_\ell = \frac{1}{2}\sum X_{\ell,ii},
$$
respectively.  In terms of the coefficients, one has:
\begin{equation} \label{both}
A_{i,jk}^\ell = \frac{1}{2(n+1)} (\delta_{ij}\delta_{lk}+\delta_{ik}\delta_{lj}) \ \ \ {\rm and}\ \ \ A_{i,jk}^\ell = \frac{1}{2} \delta_{i\ell}\delta_{jk},
\end{equation}
where $\delta$ is the Kronecker symbol.
}
\end{example}

We now show that Archimedes Lemma is automatically satisfied for any choice of coefficients $A_{i,jk}^\ell$. Let $\Delta=(V_0,\ldots,V_n)$ be a simplex and $O$ a point. Consider the simplices 
$$\Delta_i = (V_0,\ldots,V_{i-1},O,V_{i+1},\ldots,V_n),\ i=0,\ldots,n.
$$

\begin{lemma} \label{addit}
For every choice of $i,j,k$, one has:
$$
X_{i,jk}(\Delta)=\sum_{i=0}^n X_{i,jk}(\Delta_i).
$$
\end{lemma}

\proof
Let $f:\R^n \to \R^n$ be a mapping. Consider the $(n+2)\times (n+2)$-determinant
\begin{equation*} 
\begin{split}
0=\left|\begin{array}{ccccccccccccc}
f(O) & f(V_0) & f(V_1) & \cdots & f(V_n)\\
1 & 1 & 1 & \cdots & 1\\
1 & 1 & 1 & \cdots & 1
\end{array}\right| =
 \left|\begin{array}{ccccccccccccc}
f(V_0) & f(V_1) & \cdots & f(V_n)\\
1 & 1 &  \cdots & 1
\end{array}\right|  \\
- \left|\begin{array}{ccccccccccccc}
f(O) & f(V_1) & \cdots & f(V_n)\\
1 & 1 &  \cdots & 1
\end{array}\right| + \ldots + (-1)^{n+1} 
 \left|\begin{array}{ccccccccccccc}
f(O) & f(V_0) & \cdots & f(V_{n-1})\\
1 & 1 &  \cdots & 1
\end{array}\right|.
\end{split}
\end{equation*}
Let 
$$
f(x_1,\ldots,x_n)=(x_1,\ldots,x_{i-1},x_jx_k,x_{i+1},\ldots,x_n).
$$
Then, taking the orientations of the simplices into account, the above equality for determinants yields the result.
\proofend

In the next section, we shall use assumption 1, namely, the fact that the map $\Delta \mapsto C(\Delta)$ commutes with parallel translations, rotations, and transpositions of coordinates, to conclude that the coefficients $A_{i,jk}^\ell$ must be affine combinations of the ones in (\ref{both}).

\section{Proof of Theorem} \label{tedious} 

Introduce infinitesimal parallel translations in $r$th direction and rotations in the $p,q$-plane:
$$
\xi_r = \frac{\partial}{\partial x_r},\ \ \eta_{pq}=x_p \frac{\partial}{\partial x_q} - x_q \frac{\partial}{\partial x_r},\ \ p,q,r=1,\ldots,n.
$$
Let $\sigma$ denote a transposition of coordinates. The next lemma describes the action of these transformations on the polynomials $X_{i,jk}$.

\begin{lemma} \label{action}
One has
\begin{equation*}
\begin{split}
\sigma (X_{i,jk}) = - X_{\sigma(i),\sigma(j)\sigma(k)}, \
\xi_r (X_{i,jk})=(\delta_{ij}\delta_{rk}+\delta_{ik}\delta_{rj}) V,\\ 
\eta_{pq} (X_{i,jk})= \delta_{qj} X_{i,pk} + \delta_{qk} X_{i,pj} - \delta_{pi} X_{q,jk}
- \delta_{pj} X_{i,qk} - \delta_{pk} X_{i,qj} + \delta_{qi} X_{p,jk}.
\end{split}
\end{equation*}
\end{lemma}

\proof
The first equality follows from the fact that, along with the transposition of indices, exactly two rows of the determinant  $X_{i,jk}$ are interchanged. For the second equality, notice that 
$$
X_{i,jk} = x_jx_k \frac{\partial}{\partial x_i} (V),\ 
\frac{\partial}{\partial x_r} (V) =0,\ x_j \frac{\partial}{\partial x_i} (V) = \delta_{ij} V.
$$
It follows that 
$$
\xi_r (X_{i,jk})= \left[\frac{\partial}{\partial x_r},x_jx_k \frac{\partial}{\partial x_i}\right] (V) = 
\delta_{rj}\ x_k \frac{\partial}{\partial x_i} (V) + \delta_{rk}\ x_j \frac{\partial}{\partial x_i} (V)
= (\delta_{ij}\delta_{rk}+\delta_{ik}\delta_{rj}) V.
$$
The third equality is proved similarly.
\proofend

The covariance of the map $\Delta \mapsto C(\Delta)$ with respect to rigid motions is expressed by the next equations on the coefficients $A_{i,jk}^\ell$.

\begin{proposition} \label{eqcoeff}
For every transposition $\sigma$ of the indices $1,\ldots,n$, one has
\begin{equation} \label{trasp}
A_{i,jk}^\ell = A_{\sigma(i),\sigma(j)\sigma(k)}^{\sigma(\ell)}.
\end{equation}
The covariance with respect to infinitesimal translations is given by
\begin{equation} \label{transl}
\sum A_{i,ir}^\ell = \frac{1}{2} \delta_{\ell r}\ \ {\rm for\ all}\ \ \ell, r,
\end{equation}
and with respect to infinitesimal rotations by 
\begin{equation} \label{rota}
\begin{split}
A_{a,qc}^\ell \delta_{pb} + A_{a,qb}^\ell \delta_{pc} - A_{p,bc}^\ell \delta_{qa}
- A_{a,pc}^\ell \delta_{qb} - A_{a,pb}^\ell \delta_{qc} + A_{q,bc}^\ell \delta_{pa}\\
- A_{a,bc}^p \delta_{\ell q} + A_{a,bc}^q \delta_{\ell p}=0,
\end{split}
\end{equation}
for all $a,b,c,p,q,\ell$.
\end{proposition}

\proof
A transposition of coordinates is reflection in a hyperplane, and it changes the sign of $V$. Hence the covariance  of the map $\Delta \mapsto C(\Delta)$ with respect to $\sigma$ implies the equality
$$
\sum A_{i,jk}^\ell X_{\sigma(i),\sigma(j)\sigma(k)} = \sum A_{i,jk}^{\sigma(\ell)} X_{i,jk}
$$
for all $\ell$. Since the polynomials $X_{i,jk}$ form a basis, for each term on the right, there is a matching term on the left:
$$
A_{i,jk}^{\sigma(\ell)} = A_{\sigma(i),\sigma(j)\sigma(k)}^\ell.
$$
Renaming $\sigma(\ell)$ by $\ell$, we obtain (\ref{eqcoeff}).

To establish (\ref{transl}), we use Lemma \ref{action} to calculate:
$$
\xi_r \left(\frac{y_\ell}{V}\right)= \sum A_{i,jk}^\ell (\delta_{ij}\delta_{rk}+\delta_{ik}\delta_{rj})  = \sum 2 A_{i,ir}^\ell.
$$
On the other hand, $y_\ell/V$ is  the $\ell$th component of the map $\Delta \mapsto C(\Delta)$, and the infinitesimal translation in the $r$th direction sends it to $\delta_{\ell r}$.
By translation covariance,  the above sum equals $\delta_{\ell r}$, as claimed.

Likewise, the infinitesimal rotation in the $p,q$-plane annihilates $y_\ell$ for $\ell$ distinct from $p,q$, and sends $y_q$ to $y_p$, and $y_p$ to $-y_q$; in short,  
$$
y_\ell \mapsto y_p \delta_{\ell q} - y_q \delta_{\ell q}.
$$
On the other hand, by Lemma \ref{action}, 
\begin{equation*}
\begin{split}
\eta_{pq} (y_\ell) = \sum A_{i,jk}^\ell (\delta_{qj} X_{i,pk} + \delta_{qk} X_{i,pj} - \delta_{pi} X_{q,jk}
- \delta_{pj} X_{i,qk} - \delta_{pk} X_{i,qj} + \delta_{qi} X_{p,jk})\\
 = 2A_{i,qj}^\ell X_{i,pj} - A_{p,jk}^\ell X_{q,jk} - 2A_{i,pj}^\ell X_{i,qj} + A_{q,jk}^\ell X_{p,jk}.
\end{split}
\end{equation*}
Equate this to 
$$
y_p \delta_{\ell q} - y_q \delta_{\ell q} = \sum (A_{i,jk}^p \delta_{\ell q} - A_{i,jk}^q \delta_{\ell p}) X_{i,jk},
$$
and then, for fixed $a,b,c$, equate the coefficients in front of $X_{a,bc}$ in both expressions to obtain (\ref{rota}). 
\proofend

Now we need to solve the  system of linear equations (\ref{trasp})--(\ref{rota}) on the unknowns $A_{i,jk}^\ell$. We use (\ref{trasp}) to reduce the number of variables.

Consider the following four cases. If $|\{i,j,k,\ell\}| = 4$ then, applying an appropriate sequence of transpositions, we obtain: $A_{i,jk}^\ell= A^1_{2,34}=:t$. If $|\{i,j,k,\ell\}| = 3$, then one has four sub-cases, and $A_{i,jk}^\ell$ is equal to
$$
A^1_{1,23}=:u,\ {\rm or}\ A^1_{2,13}=:v,\ {\rm or}\ A^2_{3,11}=:w,\ {\rm or}\ A^2_{1,13}=:s.
$$
Likewise, if $|\{i,j,k,\ell\}| = 2$, then one has five sub-cases, and $A_{i,jk}^\ell$ is equal to
$$
A^1_{1,22}=:\phi,\ {\rm or}\ A^1_{2,12}=:\psi,\ {\rm or}\ A^2_{1,11}=:\alpha,\ {\rm or}\ A^1_{1,12}=:\beta,\ {\rm or}\ A^1_{2,11}=:\gamma.
$$
Finally, if $|\{i,j,k,\ell\}| = 1$, then $A_{i,jk}^\ell= A^1_{1,11}=:\nu$. Thus we have 11 unknowns.

Now the strategy is to consider particular cases of (\ref{rota}) and (\ref{transl}). 

To start with, consider (\ref{rota}) with $\ell=q=b=c \neq p=a$. One obtains
$$
- 2 A^q_{p,pq} + A^q_{q,qq} - A^p_{p,qq}=0,
$$
or 
$2\psi+\phi=\nu.$ Likewise, (\ref{transl}) with $\ell=r$ yields
$ (n-1) \psi + \nu = 1/2.$
It follows that 
\begin{equation} \label{soln}
\nu = \frac{1}{2} - (n-1) \psi,\ \phi = \frac{1}{2} - (n+1) \psi.
\end{equation}

We pause to check against Example \ref{twoknown}. For the center of mass,
$$
\psi = \frac{1}{2(n+1)},\ \nu=\frac{1}{n+1},
$$
and the rest of variables vanish; for the circumcenter of mass,
$\phi = \nu = 1/2,$
and the rest vanishes. In both cases, (\ref{soln}) holds.

To finish the proof of Theorem \ref{main}, we need to show that all variables, except $\phi, \psi, \nu$, vanish. We proceed in a similar fashion: (\ref{rota}) with $\ell=q=a=b=c \neq p$ yields
$$
\alpha+2\beta+\gamma=0,
$$
(\ref{rota}) with $\ell=q=c \neq p=a=b$ yields $\gamma=\alpha$. Hence $\beta=-\alpha$.
Next, (\ref{transl}) with $\ell\neq r$ yields
$$
(n-2)s+\alpha+\beta=0,
$$
hence $s=0$.

Next, (\ref{rota}) with $\ell=q=c$, but distinct from pairwise distinct $p,a,b$, yields $t=0$. It remains to eliminate $u,v,w$ and $\alpha$. Toward this,  (\ref{rota}) with $\ell=q=c\neq p=b$ and distinct from $a$ yields
$$
\gamma=v+w,\ \ {\rm hence}\ \ \alpha=v+w.
$$
Likewise,  (\ref{rota}) with $\ell=q=c\neq p=a$ and distinct from $b$, yields
$$
-s+\beta-u=0, \ \ {\rm hence}\ \ \alpha+u=0. 
$$
Next, (\ref{rota}) with $\ell=q=c=a$, but distinct from pairwise distinct $p,b$, yields
$$
u+v+s=0, \ \ {\rm hence}\ \ u+v=0,
$$
and (\ref{rota}) $\ell=q=c=b$, but distinct from pairwise distinct $p,a$, yields
$2v+w=0.$ We have obtained four linear equations on $u,v,w,\alpha$, and the only solution of this system is zero. This completes the proof.

\section{Final remarks} \label{rmks}

(i) Degenerate simplices can be safely ignored when calculating the center of mass: such a simplex has a finite centroid and zero volume, making no contribution to the total sum. Not so for the circumcenter of mass: although the volume of a nearly degenerate simplex tends to zero, its circumcenter may go to infinity, and the contribution to the  sum (\ref{totsum}) may  be non-negligible. The map $\varphi: \Delta \mapsto \Vol(\Delta)\ CC(\Delta)$, being polynomial in the coordinates of the vertices, is continuous.

\begin{figure}[hbtp]
\centering
\includegraphics[width=5in]{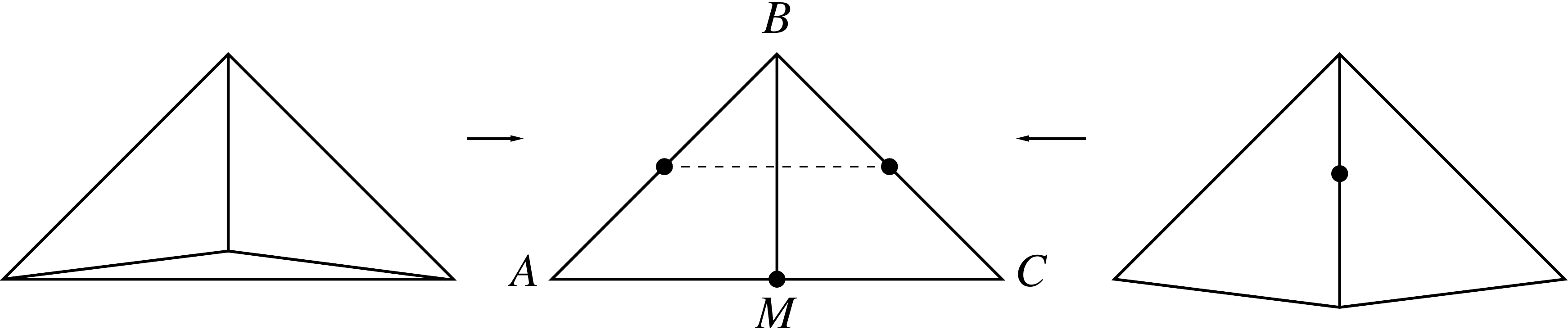}
\caption{Contribution of a degenerate triangle}
\label{discont}
\end{figure}

For example, consider an isosceles right triangle $ABC$. Its circumcenter is the midpoint $M$ of the hypothenuse $AC$. Consider the triangulation in figure \ref{discont} consisting of three triangles, one of which, $AMC$, is degenerate. If one ignored this triangle, then, by the Archimedes Lemma,  the circumcenter of mass of $\triangle ABC$ would be the midpoint of the segment connecting the midpoints of the hypothenuses $AB$ and $BC$ of the triangles $ABM$ and $BCM$. The latter point is the circumcenter of the {\it quadrilateral} $ABCM$, not the {\it triangle} $ABC$.  

\begin{figure}[hbtp]
\centering
\includegraphics[width=2in]{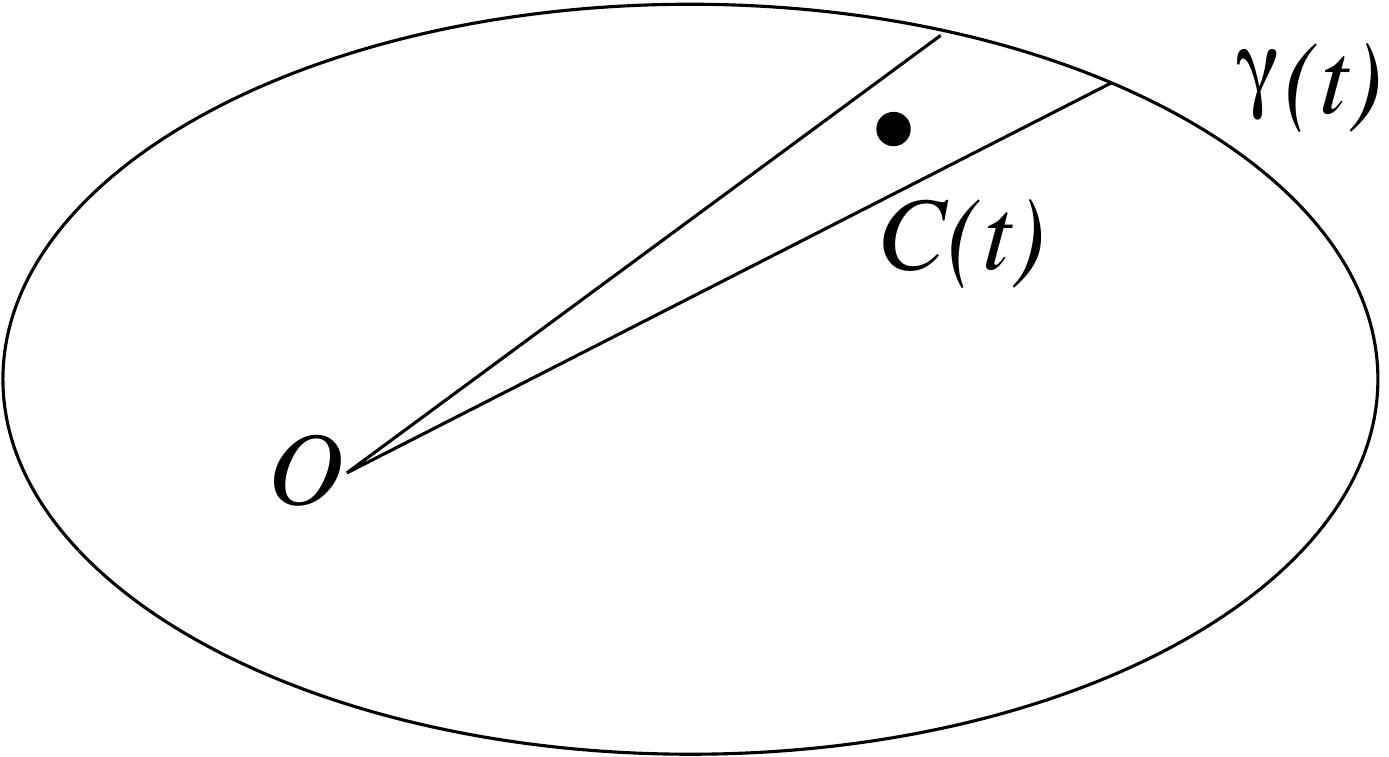}
\caption{Continuous limit of the circumcenter of mass}
\label{curve}
\end{figure}

(ii) One may wish to extend the notion of the circumcenter of mass to more general sets. For example, let $\gamma(t)$ be a parameterized smooth curve, star-shaped with respect to point $O$, see figure \ref{curve}.
It is natural to define the Circumcenter of Mass by continuity as
\begin{equation} \label{integ}
\frac{\int {C}(t)\ dA}{\int dA},
\end{equation}
where $C(t)$ denotes the limiting $\varepsilon \to 0$ position of the vector from $O$ to the circumcenter of the infinitesimal triangle $O \gamma(t) \gamma(t+\varepsilon)$, and $dA$ is the area of this infinitesimal triangle. However, this  does not give anything new: the integral (\ref{integ}) is  the center of mass of the lamina bounded by the curve \cite{TT2}.

(iii) Although the rational map $\Delta \mapsto CC(\Delta)$ is discontinuous, the polynomial map $\varphi : P \mapsto \Vol(P)\ CCM(P)$, defined on simplicial polytopes in $\R^n$, is continuous and is a valuation:
$$
\varphi (P_1 \cup P_2) + \varphi (P_1\cap P_2) = \varphi(P_1) + \varphi(P_2). 
$$
This valuation is isometry covariant; see, e.g., \cite{Sch} for the theory of valuations. 

$\R^n$-valued 
continuous isometry covariant valuations on convex compact subsets of $\R^n$ were classified in \cite{HS} as linear combinations of intrinsic moments. Namely, given a convex set $K$, let $K_\varepsilon$ be its $\varepsilon$-neighborhood. Then the moment vector
\begin{equation} \label{moment}
\int_{K_\varepsilon} x\ dx
\end{equation}
is a polynomial in $\varepsilon$, and its coefficients span the space ${\cal V}_n$ of continuous isometry covariant valuations (the free term being $\Vol(K)\ CM(K)$). One has: dim ${\cal V}_n = n+1$.

The next proposition states that the circumcenter of mass is not a linear combination of the intrinsic moments.

\begin{proposition} \label{new}
The map $\varphi : P \mapsto \Vol(P)\ CCM(P)$ is not an element of the space ${\cal V}_n$.
\end{proposition}

\proof
We argue in dimension 2; the general case is similar.

Consider an isosceles triangle $K(\alpha)$ whose base is aligned with the $x$-axis and has length 2, whose axis of symmetry is the $y$-axis, and whose base angle is $\alpha$. If $\alpha$ is close to zero then the moment vector (\ref{moment}) is close to the zero vector, and so are  the coefficients of the powers of $\varepsilon$ in (\ref{moment}) that constitute a basis in ${\cal V}_2$. 

Assume that $\varphi$ is a linear combination of the three basic vectors of the space ${\cal V}_2$ (with constant coefficients, independent of $K$). Then $\varphi(K(\alpha)) \to 0$ as $\alpha\to 0$. However, a straightforward computation shows that 
$$
\lim_{\alpha\to 0} \varphi(K(\alpha)) = \left(0, -\frac{1}{2}\right)
$$
(without computation, this can be seen in figure \ref{discont}). 
This is  a contradiction. 
\proofend

(iv) We finish with a conjecture and two problems. 

Assume that one assigns a ``center" to every simplicial polytope in $\R^n$ so that the center depends analytically on the polytope, commutes with dilations, and satisfies the Archimedes Lemma (with the weights equal to the respective volumes).

\begin{conjecture} \label{conj}
The space of such centers is 1-dimensional: they are affine combinations of the centers of mass and the circumcenters of mass.
\end{conjecture}

For $n=2$, this is proved in \cite{TT2}.

As we mentioned, the construction of the circumcenter of mass extends to the spherical and hyperbolic geometries \cite{Ak,TT2}. Thus one has versions of Conjecture \ref{conj} for $\Sph^n$ and $\HH^n$ as well.

Next, we pose the following problem.

\begin{problem} \label{prbl}
Describe $\R^n$-valued continuous isometry covariant valuations on simplicial polytopes in $\R^n$.
\end{problem}

Finally, it is interesting to find an axiomatic description of the centers for simplicial polygons and polytopes, discussed in this note, in the  three geometries of constant curvature (for the center of mass, see \cite{Ga}). These centers should be isometry covariant and satisfy some additivity condition (Archimedes Lemma or valuation-like). 
Such a description should include the valuations from Problem \ref{prbl}. At the moment of writing, we do not know such an axiomatic description.

\bigskip
{\bf Acknowledgments}. This work is an extension of the project that originated in the  Summer@ICERM 2012 program; it is a pleasure to acknowledge the inspiring atmosphere and hospitality of the institute. We are grateful to V. Adler, A. Akopyan, Yu. Baryshnikov and B. Gr\"unbaum for their interest and help.
The first author was  supported by the NSF grant DMS-1105442, and the second author was supported by  a NSF Graduate Research Fellowship under Grant No. DGE 1106400

\end{document}